\documentclass[11pt]{amsart}

\usepackage{latexsym}
\usepackage{color}

\usepackage{tikz-cd} \usepackage[all]{xy} \usepackage[utf8x]{inputenc}
\usepackage[english]{babel} \usepackage{mathtools}
\usepackage[T1]{fontenc} \usepackage{lmodern} \usepackage{mathrsfs}
\usepackage{enumitem}

\usepackage{graphicx,epsfig,amscd,graphics,xypic,subcaption}
\usepackage[all]{xy} \usepackage{amsmath,amssymb}
\usepackage[margin=4cm]{geometry}

\usepackage{hyperref} \hypersetup{
  bookmarks=true, 
  unicode=false, 
  pdftoolbar=true, 
  pdfmenubar=true, 
  pdffitwindow=false, 
  pdfstartview={FitH}, 
  pdftitle={Lyapunov exponents of the Hodge bundle over strata of
    quadratic differentials with large number of poles}, 
  pdfauthor={Charles Fougeron}, 
  pdfkeywords={Lyapunov exponents} {Quadratic differentials} {large
    number of poles}, 
  pdfnewwindow=true, 
  colorlinks=false, 
  linkcolor=red, 
  citecolor=green, 
  filecolor=magenta, 
  urlcolor=cyan 
}

\newcommand{\sheaf}[2]{\mathcal{O}_{#1} \negmedspace \left( #2
  \right)} \newcommand{\sheafx}[1]{\mathcal{O} \negmedspace
  \left({#1}\right)} \newcommand{\tsheaf}[1]{\mathcal{O}_{#1}}

\newcommand{\abs}[1]{\left| #1 \right|}
\newcommand{\QuadraticStratum}[1]{\mathcal{Q} \left( #1 \right)}
\newcommand{\PQuadraticStratum}[1]{P \mathcal{Q} \left( #1 \right)}

\newcommand{\rank}{\operatorname{rk}} \newcommand{\ff}{f_* \,}
\newcommand{\fb}{f^*}

\newcommand{\zandp}{\sum_{i=1}^n m_i Z_i - \sum_{j=1}^p P_i}
\newcommand{\onlyz}{\sum_{i=1}^n m_i Z_i}

\newcommand{\compactX}{\overline {\mathcal X}}
\newcommand{\X}{\mathcal X}

\newcommand{\dualizing}{\omega_{\compactX / \overline B}}

\newcommand{\C}{\mathbb{C}} 
\newcommand{\Q}{\mathbb{Q}} \newcommand{\R}{\mathbb{R}}
\newcommand{\GL}{\mathrm{GL}} \newcommand{\SL}{\mathrm{SL}}
\newcommand{\SO}{\mathrm{SO}} \newcommand{\Id}{\operatorname{Id}}

\newtheorem*{theorem*}{Theorem} 
\newtheorem*{corollary*}{Corollary}
 
\newtheorem{lemma}{Lemma}

\newtheorem*{conjecture*}{Conjecture}

\theoremstyle{definition}

\newtheorem*{remark}{Remark} 

\begin{document}

\title[Lyapunov exponents with large number of poles] {Lyapunov
  exponents of the Hodge bundle over strata of quadratic differentials
  with large number of poles}

\author{Charles Fougeron}

\maketitle

\begin{abstract}
  We show an upper bound for the sum of positive Lyapunov exponents of
  any Teichmüller curve in strata of quadratic differentials with at
  least one zero of large multiplicity. As a corollary it holds for
  any $\SL(2,\R)$-invariant submanifold defined over $\Q$ in these strata.
  This proves Grivaux-Hubert's
  conjecture about the asymptotics of Lyapunov exponents for strata
  with a large number of poles in the situation when at least one zero
  has high multiplicity.
\end{abstract}

\section{Introduction.}
Lyapunov exponents of translation surfaces were introduced two
decades ago by Zorich in \cite{zo1} and \cite{zo2} as dynamical
invariants which describe how associated leaves \textit{wind around
  the surface}.  On any translation surface we can introduce a
translation flow which generalizes the linear flow on a flat torus
(see \cite{zo3} for an introduction to the subject). This
flow has a very simple local dynamic --- it is a parabolic system.
Nonetheless the homology of the translation flow presents a rich
asymptotic behaviour and its deviation from the asymptotic cycle is
described by Lyapunov exponents of the Hodge bundle over a
$\SL_2(\R)$ invariant subspace in the moduli space of curves.

Even if a numerical approximation of these exponents is accessible
\cite{flatsurf}, there is \textit{a priori} no hope for an explicit formula to
compute them.  Yet a breakthrough of Eskin, Kontsevich and Zorich
showed an astonishing formula binding the sum of the Lyapunov
exponents to the Siegel-Veech constant of the invariant locus
\cite{EKZ}. This theorem followed an insightful observation of \cite{K} that
this sum is related to the degree of a Hodge subbundle, which
was proven later in \cite{Forni}.  This was a starting point to
evaluate Lyapunov exponents in certain particularly symmetric cases,
for example for square-tiled cyclic covers \cite{FMZ}, \cite{EKZ2} and
triangle groups \cite{BouwMoeller}; and to compute explicitly diffusion
rate of wind-tree models \cite{DHL}, \cite{DZ}.

Other advances have since been made to estimate Lyapunov exponents for
higher genus.  In \cite{fei3}, Yu gave a partial proof of the conjecture
of \cite{KZ} that the second Lyapunov exponent for hyperelliptic
components of strata should go to 1 as the genus goes to infinity. His
proof was conditional on a conjecture he introduced in the same article
which brought new ideas to find
bounds from below to these exponents. This conjecture has recently been
proven in \cite{EKMZ} (see \cite{DD} for a foliation-theoretic
point of view).  Yu also obtained an upper bound for the sum of Lyapunov
exponents with respect to Weierstrass gaps.  Yu's idea exposed in
\cite{fei1} and \cite{fei2} that independently appeared in \cite{CM1}
and \cite{CM2} is to use algebraic characterization of sum of Lyapunov
exponents and estimate it with homological algebra arguments.

In parallel Grivaux and Hubert remarked that some Teichmüller curves
whose Lyapunov exponents are all zero can appear in quadratic
strata. Moreover in \cite{GH} they prove that for this to happen the
curve should be in a stratum with at least $\max(2g-2,2)$ poles.  The
heuristic explanation they provide for this result is that
\textit{poles slow down the linear flow}, since passing by a pole 
makes the flow to retrace its steps. This leads them to the
following conjecture transmitted to the author by oral communication.

\begin{conjecture*}[Grivaux-Hubert]
  Positive Lyapunov exponents associated to a translation surface have
  a uniform bound depending only on the number of poles of the
  surface. This bound goes to zero when the number of poles goes to
  infinity.
\end{conjecture*}

Pascal Hubert pointed out later on that Remark 2 of \cite{DZ} gives a
counter-example to the conjecture in this generality.
Taking a cover of a base surface, Delecroix and Zorich exhibited
a family of surfaces of genus 1 with an arbitrarily large number of poles,
in which the first Lyapunov exponent is always equal to $\frac 2 3$. 

Nevertheless, in \cite{CM2} Theorem 8.1, D. Chen and M. Möller obtained
a result in this direction for $\QuadraticStratum{n,-1^n}$ and
$\QuadraticStratum{n,1,-1^{n+1}}$, showing that these strata are
non-varying and computing explicitly the sum of Lyapunov exponents for
every Teichmüller curve which is equal to $2/(n+2)$.

Inspired by Yu's homological methods we obtain
a general upper bound in terms of the highest multiplicity of
zeros and genus. This proves the conjecture in this case.

\begin{theorem*}
  \label{main}
  For any Teichmüller curve $\mathcal C$ in a quadratic stratum
  $\QuadraticStratum{m_1, \dots, m_k, -1^p}$ of genus $g$ where
  $m_1 \geq m_2 \geq \dots \geq m_k$, if $m_1 \geq 2g$, then
  $$ L^+(\mathcal C) \leq \frac {(3g-1)g}{m_1+2}, $$
  where $L^+(\mathcal C)$ stands for the sum of its positive Lyapunov
  exponents.
\end{theorem*}

Using recent advances in the theory of $\SL(2,\R)$-invariant
submanifold, we obtain as a corollary a very general statement.

\begin{corollary*}
  \label{coro}
  For any half-translation surface in the stratum
  $\QuadraticStratum{m_1, \dots, m_k, -1^p}$ of genus $g$ where
  $m_1 \geq m_2 \geq \dots \geq m_k$, if $m_1 \geq 2g$ and its
  $\SL(2,\R)$-orbit closure contains some square tiled surface, then
  the following inequality holds:
  $$ L^+(\mathcal N) \leq \frac {(3g-1)g}{m_1+2}, $$
  where $L^+(\mathcal C)$ stands for the sum of the positive Lyapunov
  exponents associated to the surface.
\end{corollary*}

Another direction in which the conjecture can be studied is in the
case of whole strata.  Developing a tool to compute Lyapunov
exponents for quadratic differentials (available in \cite{flatsurf}),
we have performed broad computer experiments on strata with number of
zeros and poles both going to infinity. These experiments tend to show that the sum
of their Lyapunov exponents goes to zero at speed $\frac{1}{\sqrt{p}}$.
Consequently, we make the following conjecture,

\begin{conjecture*}
  Let $Q_n$ be a sequence of connected components of strata for a
  fixed genus and $p_n$ their number of poles. If $p_n$ 
  goes to infinity, then,
  $$ \lambda^+_1(Q_n) = \mathcal O (1/\sqrt{p_n}). $$
\end{conjecture*}

\subsection*{Acknowledgements.}
I am very grateful to Fei Yu for very enlightening and friendly
discussions in MFO and Beijing. I thank Pascal Hubert and Julien
Grivaux who took time to discuss their conjecture and this work in a
most pleasant atmosphere in Marseille.

\section{Background material.}
\subsection*{Strata.}
A half-translation surface is a pair $(X,q)$ where $X$ is a Riemann
surface and $q$ is a quadratic differential on $X$ with at worst 
simple poles.  If $S(q)$ is the set of zeros and poles of $q$ on $X$,
we can endow $\tilde X := X \backslash S(q)$ with charts
$\phi_i: U_i \rightarrow X$ such that $\phi_i^* q = dz^{\otimes 2}$.
In such an atlas, the transition functions are translations
composed with $\pm \Id$.  The quadratic differential induces a flat
metric $\abs{q}$ on the open surface $\tilde X$. This metric can be
extended to the whole surface and has conical points at $S(q)$
with angles multiples of $\pi$. If we fix integers
$m_1, \dots, m_k, p$ such that $\sum_{i=1}^k m_i - p = 4g-4$ for some
positive integer $g$, the stratum of half-translation surfaces
$\QuadraticStratum{m_1, \dots, m_k, -1^p}$ is the set of
half-translation surfaces $(X, q)$ where $q$ has $k$ distinct zeros with
multiplicity $m_1, \dots, m_k$ and $p$ simple poles.  The
projectivized space $\PQuadraticStratum{m_1, \dots, m_k, -1^p}$ is
obtained by taking the quotient under the scalar action of $\C^*$
on the quadratic differential.\\



\subsection*{Teichmüller curves.}
There is a natural action of $\GL(2,\R)$ on each stratum. In the
flat atlas above, the charts were constructed in such a way that the
transition functions are translations by vectors $v_{ij}$ composed
with $\pm \Id$.  For a matrix $M \in \GL(2,\R)$ we define a new
surface by multiplying the previous change of charts by $M$. These new
transition maps will be translations by vectors $Mv_{ij}$ composed with $\pm \Id$.

Let $\mathbb H = \SL(2,\R)/\SO(2,\R)$ denote the Poincaré upper
half-plane.  For any $(X,q)$ the action of $SL(2,\R)$ factors to a map
$$\mathbb H \rightarrow \PQuadraticStratum{m_1, \dots, m_k, -1^p}$$
which is an immersion. The image of this map is called a Teichmüller
disk.
This map also factors through the Veech group which is its stabilizer.
If the Veech group $\Gamma$ of a given half-translation surface is a
lattice in $\SL(2,\R)$ we say that the surface is a Veech surface and
the image of $\mathbb H/\Gamma$ in the projective stratum is called
a Teichmüller curve.  Another point of view for the Teichmüller curve
is to consider the induced algebraic map
$\mathcal C \rightarrow \mathcal M_g$. This is the convention used in
\cite{Moller}.

By definition, a Teichmüller curve is a surface with hyperbolic
structure --- which may have cusps --- parametrizing a family of curves. 
A main tool for studying these
curves is their compactification at the cusps.  If we consider such a
family $f: \X \rightarrow C$, the stable reduction theorem (see
\cite{HM} Theorem 3.47) asserts that there exists a finite-index
cover $\overline B \to \overline C$ ramified over the boundary of $C$
such that we can define a family
$\overline f: \compactX \rightarrow \overline B$ which extends $f$,
and such that the fibers at the boundary of $B$ are stable curves.

The structure sheaf of $\compactX$ will be denoted by $\mathcal O$,
the dualizing sheaf of
$\overline {\mathcal X} \rightarrow \overline B$ by $\dualizing$, and
the opposite of the Euler characteristic of $\overline B$ by $\chi$.


\subsection*{Divisors classes}
For a Teichmüller curve generated by a quadratic differential $q$,
there is an associated line bundle
$\mathcal L \subset \ff \dualizing^{\otimes2}$ on $\overline B$ whose
fiber over the point corresponding to $(X,q)$ is $\C \cdot q$.  This
bundle has a maximality property (see \cite{Moller} section 2, or Theorem of \cite{Goujard}) 
which implies that $\deg(\mathcal L) = \chi$.


We denote by $Z_1, \dots, Z_k$ and $P_1, \dots, P_p$ sections of
$\compactX \rightarrow \overline B$, which intersect the fibers at
each zero with multiplicity $m_1, \dots, m_k$ and at poles with multiplicity $1$.  For
convenience we also introduce the total divisor $\mathcal A := \zandp$
and the divisor of zeros $\mathcal Z := \onlyz$.

Proposition 4.2 of \cite{CM2} gives a formula for the self-intersection of $Z_i$. 
We reproduce the argument of this computation. First notice that
$\mathcal L \subset f_* \dualizing^{\otimes 2} (P_1+ \dots+ P_p)$.
The pullback of this inclusion yields
$$0 \rightarrow \fb \mathcal L \rightarrow \dualizing^{\otimes 2}(P_1+ \dots + P_p) \rightarrow \sheaf{\mathcal Z}{\mathcal Z} \rightarrow 0.$$
This is true for non-compactified Teichmüller
spaces, and it remains true after compactification since the
multiplicities of the zeros of the limit differentials stay the same
on stable curves (see section 5.4 of \cite{CM2} for details on
cusps of Teichmüller curves). This yields an isomorphism
$\fb \mathcal L (\mathcal A) \simeq \dualizing^{\otimes 2}.$\\

For completeness, we reproduce below the computation of
Proposition 4.2 in \cite{CM2}.
According to the adjunction formula, for any singularity $Z_i$,
$$ \dualizing(Z_i)_{|Z_i} = \omega_{Z_i} = 0,$$ which implies
$Z_i^2 = - \dualizing \cdot Z_i$.  Moreover $Z_i \cdot Z_j = 0$ for
any $i \neq j$, since two distinct fibers of $f$ do not intersection.
Finally,
$$2 Z_i^2 = - m_i Z_i^2 - \deg \mathcal L,$$
where $m_i$ is the corresponding multiplicity.\\

\noindent Which leads us to the formula:
\begin{equation}
  \label{autointersection}
  Z_i^2 = \frac {-\chi}{m_i+2}.
\end{equation}



\section{Lyapunov exponents}
In the remainder of the paper, $\mathcal C$
will stand for a Teichmüller curve in
$$\QuadraticStratum{m_1, \dots, m_k, -1^p}.$$
Lyapunov exponents are dynamical invariants for translation surfaces
which describe the behaviour of the translation flow on them. For an
introduction to their diverse aspects in translation surfaces see
\cite{zo3}. We will adopt an algebraic point of view to compute
them. Our strategy is based on the observation of \cite{K} that their sum is
related to the degree of a Hodge subbundle on a Teichmüller curve (see
\cite{EKZ} (3.11) or \cite{Forni} for a proof).
We will use in the following a more precise formula computed in
\cite{CM2} as (20),
$$
L^+(\mathcal C) = \frac 2 \chi \cdot \deg \ff \sheafx{\mathcal A} +
(6g-6) - \frac 1 2 \left(\sum_{i=1}^n \frac {m_i(m_i+4)} {m_i+2} -
  3p\right),
$$
where $L^+(\mathcal C)$ is the sum of the positive Lyapunov exponents
associated to the Teichmüller curve $\mathcal C$. Recall that
Lyapunov exponents can be defined on both bases $C$ and $B$ and that
they will agree. Moreover, here the degrees are always
computed on vector bundles over the compact curve $\overline B$.\\


By Gauss-Bonnet,
$$\sum_{i=1}^n m_i - p = 4g - 4,$$
$$6g-6 + \frac 3 2 p = \frac 3 2 \sum_{i=1}^n m_i.$$
Notice that $m_i(m_i+4) - 3m_i(m_i+2) = -2m_i(m_i+1)$, leading to
the following formula which will be convenient for us:
\begin{equation}
  \label{sum}
  L^+(\mathcal C) = \frac 2 \chi \cdot \deg \ff \sheafx{\mathcal A} 
  + \sum_{i=1}^n \frac {m_i(m_i+1)} {m_i+2}.
\end{equation}
This equation relates the sum of Lyapunov exponents for a Teichmüller
curve to an invariant of the stratum and the Chern class of a sheaf,
describing the behaviour of the curve at its cusps.

To estimate the sum of the Lyapunov exponents, we are reduced to bounding the
degree of the sheaf $\ff \sheafx{\mathcal A}$. This will be done based
on the two following lemmas, 

\begin{lemma}
  $$\deg \ff \sheafx{\mathcal A} = \deg \ff \sheafx{\mathcal Z}.$$
\end{lemma}

\begin{proof}
The main tool in the introduction of this filtration is the structural
short exact sequence associated to any divisor $D$ in $\compactX$
$$ 0 \rightarrow \sheafx{-D} \rightarrow \mathcal O \rightarrow \tsheaf{D} \rightarrow 0 $$
Assume here that $D$ is one of the divisors $P_j$.  We tensor this
exact sequence by $\fb \mathcal L \otimes \sheafx{\mathcal A + D}$,
where $\mathcal A$ is the total divisor as introduced in section 2.
The functor $\ff$ gives a long exact sequence,
$$
\begin{array}{lll}
  0 &\rightarrow& \ff (\omega^{\otimes 2})  
  \rightarrow \mathcal L \otimes \ff \sheafx{\mathcal A + D}
  \rightarrow \mathcal L \\
  &\overset \delta \rightarrow& R^1 \ff (\omega^{\otimes 2}) 
  \ \longrightarrow \ R^1 \ff \left( \fb \mathcal L \otimes \sheafx{\mathcal A + D} \right)
  \ \longrightarrow \ 0
\end{array}
$$
Where we use the fact that $f$ induces an isomorphism between $D$ and
$\compactX$, which implies the
equalities $\ff \tsheaf{D} = \tsheaf{\compactX}$ and $\omega^{\otimes 2} = \fb \mathcal L \otimes \sheafx{ \mathcal A}$.\\
				      
Note that $R^1 \ff (\omega^{\otimes 2})$ over each stable curve $X$ of the
compactified Teichmüller curve has dimension
$$h^1(X, \omega_X^{\otimes 2}) = h^0(X, \omega^{\otimes -1}_X) = 0$$ by 
definition of the dualizing sheaf (see \cite{Hartshorne} III.7).
Hence this sheaf is zero, and so is
$R^1 \ff (\fb \mathcal L \otimes \sheafx{\mathcal A + D})$.  We end up
with a short exact sequence which implies a relation on the degrees,
$$\deg \left(\mathcal L \otimes \ff \sheafx{\mathcal A + D} \right) = 
\deg \left( \mathcal L \otimes \ff \sheafx{\mathcal A} \right) + \deg
\left( \mathcal L\right).$$
We reiterate the previous reasoning with $\fb \mathcal L \otimes \sheafx{\mathcal A + D}$,
as we noticed its image under the derived functor $R^1 \ff$ is zero
and so will be the image of
$\fb \mathcal L \otimes \sheafx{\mathcal A + D + D'}$ for any other
divisor $D'$ picked amoung the $P_j$.  We do so until all the pole
divisors are eliminated, and get a degree formula about the
divisor without poles $\mathcal Z$,
\begin{equation*}
	\deg \left(\mathcal L \otimes \ff \sheafx{\mathcal Z} \right) = 
	\deg \left( \mathcal L \otimes \ff \sheafx{\mathcal A} \right) + p \cdot \deg \left( \mathcal L\right).
\end{equation*}
To finish, observe that the fibers of $\ff \sheafx{\mathcal Z}$ have
constant dimension $3g-3+p$ by Riemann-Roch theorem, since the degree
of the divisor on any fiber will be $4g-4 + p > 2g-2$ for
$g, p \geq 1$.  By the Grauert Semicontinuity Theorem (see 
\cite{Hartshorne} III.12.9) this sheaf will be a vector bundle, and thus, so will be
$\ff \sheafx{\mathcal A}$, as it is a subsheaf of the latter. Moreover
$\rank \ff \sheafx{\mathcal Z} = 3g-3+p$ and $\rank \ff \sheafx{\mathcal A} = 3g-3$.\\

The previous formula becomes
$$
(3g-3+p) \deg \left(\mathcal L \right) + \deg \left(\ff
\sheafx{\mathcal Z} \right) = (3g-3) \deg \left( \mathcal L \right)
+ \deg \left( \ff \sheafx{\mathcal A} \right) + p \cdot \deg \left(
\mathcal L\right)
$$
And finally,
$\deg \left(\ff \sheafx{\mathcal Z} \right) =  \deg \left(\ff \sheafx{\mathcal A} \right).$
\end{proof}

Thus the pole divisors will not interfere in the degree we need to
estimate.  By a similar method, we will be able to bound the degree
while taking out zeros one after the other. We will do so until no
zero is left and end up with the structure sheaf $\mathcal O$. 

Considering divisors $D$ and $D'$ which are linear combinations of
$Z_i$'s with non-negative coefficients such that $D-D'=Z_i$, the dimension of global sections
for their corresponding sheaf generally jumps by one when passing from
$D$ to $D'$. Sometimes it remains constant, and this phenomenon is
usually referred to as Noether or Weierstrass gaps. Using the Riemann-Roch theorem
it can be shown that there are at most $g$ gaps in an infinite
sequence of such divisors (see \cite{FK} III.5.4).

\begin{lemma}
Let $D$ and $D'$ be linear combinations of $Z_i$'s with non-negative
coefficients, such that $D-D'=Z_i$ for some $i$; then,
$$ \deg \ff \sheafx{D} \leq \deg \ff \sheafx{D'} - \mu_i(D) \frac{\chi}{m_i+2} + R(D,D'),$$
where $R(D,D')$ is zero if there is no gap and $\mu_i(D) \frac{\chi}{m_i+2}$ otherwise. 
\end{lemma}

\begin{proof}
Similarly to the previous proof, we tensor the structural exact
sequence for divisor $Z_i$ by
$\sheafx{D}$, and we get the following
long exact sequence:
$$ 0 \rightarrow \ff \sheafx{D} \rightarrow \ff \sheafx {D'} 
\rightarrow \sheaf{Z_i}{D}
\overset{\delta}{\rightarrow} \dots$$ This gives the formula,
$$ \deg \ff \sheafx{D} - \deg \ff \sheafx{D'} = \deg \ker \delta.$$
As $\ker \delta$ is a subsheaf of
$\sheaf{Z_i}{m_i Z_i}$ which is locally free on a
complex curve, both are locally free and when $\ker \delta$ is not
zero, \textit{i.e.}, when $i \in G(\mathcal F)^c$,
$$\deg \ker \delta \leq \deg \sheaf{Z_i}{m_i Z_i} = m_i Z_i \cdot Z_i,$$
Hence by formula (\ref{autointersection}),
$$ \deg \ff \sheafx{D} \leq \deg \ff \sheafx{D'} 
- m_i \frac{\chi}{\mu^i+2}.$$
If $i \in G(\mathcal F)$, $\ker \delta = 0$, and thus
$$\ff \sheafx{D} \simeq \ff \sheafx{D'}.$$
\end{proof}

If we pick a sequence of such divisors $\mathcal F = (D_j)_{1 \leq j \leq 4g-4+p}$
that removes the zeros one by one,
using Lemma 1 and formula (\ref{sum}),
\begin{equation*}
  \label{bound}
  L^+(\mathcal C) \leq \sum_{l \in G(\mathcal F)} R(D_l,D_{l+1}), 
\end{equation*}
where $G(\mathcal F)$ is the set of indices of gaps in the sequence
$\mathcal F$.

This will lead to a proof of the Theorem written in the
introduction.  Assume we are in the setting of the theorem, we denote
the larger zero multiplicity by $m_1 \geq 2g-1$. We pick a family
$\mathcal F$ such that we first remove all zeros different from the
first one in an arbitrary way, and finish by taking off all the
multiplicity of the principal one.  
Thus,
$$L^+(\mathcal C) \leq \sum_{l \in G(\mathcal F)} R(D_l,D_{l+1}) \leq \sum_{i=1}^g \frac {2(2g - i)}{m_1 + 2} 
\leq \frac {(3g-1)g}{m_1 +2},$$
where we use the fact that Weierstrass gaps happen 
for multiplicity smaller than $2g$.

\begin{remark}
  If we can show that one of the zeros with multiplicity larger than
  $g$ is not a Weierstrass point generically, the bound becomes
  $\frac {(g+1)g} {m_1 + 2}$.
\end{remark}

We finish by proving the corollary claimed in the introduction.
An \textit{affine invariant manifold} is an immersed manifold in a
stratum of abelian differentials whose image is locally defined by real
linear equations in \textit{period coordinates} (see \cite{zo3},
Chapter 3).  Each affine invariant manifold $\mathcal N$ is the
support of an ergodic $\SL(2,\R)$-invariant probability measure which
is locally in period coordinates the restriction of Lebesgue measure
on the affine subspace (see \cite{EM}).

\begin{theorem*}[\cite{BEW}, Theorem 2.8 and \cite{EMM}, Theorem 2.3]
  Let $\mathcal N_n$ be a sequence of affine invariant manifolds and
  suppose the sequence of affine measures
  $\nu_{\mathcal N_n}$ converges to $\nu$ in the weak-* topology. Then
  $\nu$ is a probability measure, and it is an affine measure
  $\nu_{\mathcal N}$ where $\mathcal N$ is the smallest submanifold
  such that there exists $n_0$ with
  $\mathcal N_n \subset \mathcal N$ for all $n>n_0$.  Moreover the
  Lyapunov exponents of $\nu_{\mathcal N_n}$ converge to the Lyapunov
  exponents of $\nu$.
\end{theorem*}

According to \cite{Wright2}[Theorem 1.1], if $\mathcal N$ contains one
square tiled surface, its field of definition is $\Q$; thus it contains
an infinite number of square tiled surfaces,
which implies corollary \ref{coro}.

\bibliographystyle{halpha} \bibliography{HG}
 
\end{document}